\documentstyle[12pt]{article}
\def\Bbb R{{\rm \bf R}}
\def\proclaim#1{\vskip2mm{\bf #1}\em}
\def\endproclaim{\em \vskip2mm}
\def\tag#1{\eqno(#1)}
\def\gathered{\begin{array}{c}}
\def\endgathered{\end{array}}
\def\text{\mbox}

\begin{document}

\title {The enclosure method for the heat equation using time-reversal invariance for a wave equation}
\author{Masaru IKEHATA\footnote{
Laboratory of Mathematics,
Graduate School of Engineering,
Hiroshima University,
Higashihiroshima 739-8527, JAPAN}}
\maketitle

\begin{abstract}
The heat equation does not have {\it time-reversal invariance}.  However, using a solution
of an associated wave equation which has time-reversal invariance, one can establish an explicit
extraction formula of the {\it minimum sphere} that is centered at an arbitrary given point and encloses an unknown
cavity inside a heat conductive body.  
The data employed in the formula consist of a special heat flux depending on 
a large parameter prescribed on the surface of the body over an
arbitrary fixed finite time interval and the corresponding temperature field.
The heat flux never blows up as the parameter tends to infinity. This is different from a previous formula for the heat equation
which also yields the minimum sphere.  In this sense, the prescribed heat flux is {\it moderate}.

\noindent
AMS: 35R30, 35L05, 35K05

\end{abstract}


\section{Introduction}
This paper is concerned with the {\it methodology} for so-called inverse obstacle problems
which aims at reconstructing or extracting an unknown discontinuity such as an 
obstacle, inclusion, cavity, etc. embedded or occurred in a reference medium for various observation data.
As one of direct analytical methods which have an exact mathematical base the author
has introduced the {\it enclosure method} in \cite{E00} and \cite{E001} using infinitely many data and a single set of data, respectively.  
The enclosure method enables us to find a domain that encloses the unknown discontinuity.

Originally the method treats a signal which does not depend on a time variable.
Clearly, this restricts the possible application of the enclosure method since
there are a lot of inverse problems which employs various data depending on a time variable.
However, in \cite{I4} the author found an idea how to treat a signal which depends on a time variable 
in a {\it bounded interval}.
Now we have various realizations of this {\it time domain enclosure method}
in three-space dimensions for the signals governed by
the heat and wave equations \cite{IW00, IFR, IEO2, IEO3, Iwall, IK1, IK2, IK3}, 
a system of equations in a linearized viscoelastic body \cite{II},
the Maxwell system \cite{IMax, IMax2} and references therein.

Quite recently, in \cite{EV} the author considered an inverse obstacle problem for the wave
governed by the wave equation in a bounded domain and introduced an idea which combines
the {\it time-reversal invariance} of the wave equation with the newest version of enclosure method \cite{EIV} 
which employs a {\it single point} on the graph of the {\it response operator} associated with the inverse obstacle problem.
Note that the idea itself in \cite{EIV} has been applied also to a system of equations in the linear theory of thermoelasticity in \cite{Ithermo}.

The key point of the idea introduced in \cite{EV} is:  making a {\it time-reversal operation} on the solution
of the free space wave equation supported on an arbitrary given ball at the initial time.
This gives us an explicit input Neumann data depending on the ball such that
the corresponding Dirichlet data as the output over a finite time interval yields the minimum sphere having the same center point as the ball
and enclosing the obstacle.

It is curious to know whether there is an analogous approach to inverse obstacle problems
for the signal governed by the heat equation.  However, one can easily understand that the heat equation does not have time-reversal invariance
and thus everyone may think that there is no hope to establish any corresponding formula.

The aim of this paper is to show that the hope can be realized {\it partially} with the help
of time-reversal invariance of an {\it associated wave equation with a growing speed}.
One can establish an explicit
extraction formula of the minimum sphere that is centered at an arbitrary given point and encloses an unknown
cavity inside a heat conductive body.

The data employed in the formula consists of a special heat flux depending on 
a large parameter prescribed on the surface of the body over an
arbitrary fixed finite time interval and the corresponding temperature field.
The heat flux is given by taking the time-reversal operation of the normal derivative of a solution of the associated wave equation.
This idea, that is, a combination of the enclosure method and prescribing a heat flux coming from a wave equation
was initiated in \cite{PW}.

Remarkably enough, the flux never blows up as the parameter tends to infinity. 
This is different from a previous formula established
in \cite{IK2} for the heat equation with a discontinuity coefficient, which also yields the minimum sphere.
The flux used there depends on a large parameter, this is a common point, however, blows up on the whole surface of the body
as the parameter tends to infinity.
In this sense, the prescribed heat flux in this paper is {\it moderate}.
In fairness, however, it should be pointed out that the formula using the special heat flux in this paper is not a complete answer to the original question itself since the flux depends on a parameter.  This is the meaning of the partial realization of the hope mentioned above.

Note that the time-reversal operation appears in the BC-method \cite{BEL, B2}, related works \cite{BKLS, DKL, O}
and a numerical approach \cite{dBK}.  However, these are concerned with inverse problems for the wave equation.
To the author's best knowledge, no one considers any application of time-reversal invariance to 
inverse obstacle problem for the heat equation.

Now let us describe the inverse obstacle problem to be considered.
We consider the problem for the signal governed by the heat equation formulated as
below. Let $\Omega$ be a bounded domain of $\Bbb R^3$ with $C^2$-boundary.
Let $D$ be a nonempty bounded open set of $\Bbb R^3$ with $C^2$-boundary such that
$\overline D\subset\Omega$ and $\Omega\setminus\overline D$ is connected.
The set $D$ is our target which is a mathematical model of an unknown cavity in a heat conductive body $\Omega$.

Given an arbitrary positive number $T$ and $f=f(x,t)$, $(x,t)\in\partial\Omega\times\,[0,\,T]$,
let $u=u_{f}(x,t)$, $(x,t)\in(\Omega\setminus\overline D)\times\,[0,\,T]$ denote the solution of the following initial
value problem for the heat equation
$$\left\{
\begin{array}{ll}
\displaystyle
(\partial_t-\Delta)u=0 & \text{in $(\Omega\setminus\overline D)\times\,]0,\,T[$,}
\\
\\
\displaystyle
u(x,0)=0 & \text{in $\Omega\setminus\overline D$,}\\
\\
\displaystyle
\frac{\partial u}{\partial\nu}=0 & \text{on $\partial D\times\,]0,\,T[$,}\\
\\
\displaystyle
\frac{\partial u}{\partial\nu}=f(x,t) & \text{on $\partial\Omega\times\,]0,\,T[$.}
\end{array}
\right.
\tag {1.1}
$$
The $f$, $u$ and $\partial_t u$ belong to $L^2(0,T;H^{-1/2}(\partial\Omega))$, $L^2(0,\,T;H^1(\Omega\setminus\overline D))$ 
and $L^2(0,\,T;H^1(\Omega\setminus\overline D)')$, respectively.
The $u$ satisfies, for a positive constant $C_T$ beging independent of $f$
$$\displaystyle
\Vert u\Vert_{L^2(0,T;H^1(\Omega\setminus\overline D))}
+\Vert\partial_tu\Vert_{L^2(0,T;H^1(\Omega\setminus\overline D)')}
+\Vert u(\,\cdot\,,T)\Vert_{L^2(\Omega\setminus\overline D)}
\le C_T\Vert f\Vert_{L^2(0,T;H^{-1/2}(\partial\Omega))}.
\tag {1.2}
$$
See \cite{DL} for more information.

The inverse obstacle problem to be considered here is: find an extraction formula
of information about the geometry of $D$ from $u_f(x,t)$ for $(x,t)\in\partial\Omega\times\,]0,\,T[$
corresponding to finitely or infinitely many $f$.

Let $B$ be an open ball centered at an arbitrary point $p\in\Bbb R^3$ with radius $\eta$.
The aim of this paper is to show that: there exists a suitable Neumann data $f$
depending on $B$ and a large parameter $\tau$ such that

$\bullet$  $f=O(1)$ in $L^2(0,T; H^{-1/2}(\partial\Omega))$ as $\tau\rightarrow\infty$;

$\bullet$  the $u_f$ on $\partial\Omega\times\,]0,\,T[$ with an aribitrary fixed $T$
as $\tau\rightarrow\infty$ yields $R_D(p)$, where
$$\displaystyle
R_D(p)=\sup_{x\in D}\,\vert x-p\vert.
$$
Note that the quantity $R_D(p)$ gives the radius of the minimum sphere centered at $p$
that encloses the unknown cavity $D$.

From (1.2) we have, for such $f$ as $\tau\rightarrow\infty$
$$\displaystyle
\Vert u_f\Vert_{L^2(0,T;H^1(\Omega\setminus\overline D))}
+\Vert\partial_tu_f\Vert_{L^2(0,T;H^1(\Omega\setminus\overline D)')}
+\Vert u_f(\,\cdot\,,T)\Vert_{L^2(\Omega\setminus\overline D)}
=O(1).
\tag {1.3}
$$
In this sense, the temperature field $u_f$ never blows up inside the body
as $\tau\rightarrow\infty$.

Note that from a combination of \cite{IK1} and \cite{IK2}, one can easily show that: 
if the heat flux
$f$ takes the form
$$\begin{array}{lll}
\displaystyle
f(x,t)=\varphi(t)\frac{\partial}{\partial\nu}v_{\tau}(x;p), & x\in\partial\Omega, & 0<t<T,
\end{array}
$$
where
$$\displaystyle
v_{\tau}(x;p)=\left\{
\begin{array}{ll}
\displaystyle
\frac{\sinh\sqrt{\tau}\,\vert x-p\vert}{\vert x-p\vert} & \text{ if $x\in\Bbb R^3\setminus\{p\}$,}\\
\\
\displaystyle
\sqrt{\tau} & \text{if $x=p$,}
\end{array}
\right.
$$
and $\varphi\in L^2(0,\,T)$ satisfies, say $\varphi(t)\sim t^m$ with an integer $m\ge 0$ as $t\downarrow 0$, then
one can extract $R_D(p)$ from $u_f$ on $\partial\Omega$ over time interval $]0,\,T[$ as $\tau\rightarrow\infty$.
Note that the function $v_{\tau}(x;p)$ of $x\in\Bbb R^3$ is an entire solution of the modified Helmholtz equation
$(\Delta-\tau)v=0$.  However, the $f$ given above blows up as $\tau\rightarrow\infty$ in the order of exponential
everywhere on $\partial\Omega$ if $p\in\Omega$, which is the most interesting case.  Thus, this choice does not satisfy
the requirement of moderateness of the state as $\tau\rightarrow\infty$.

\section{Statement of the result and its proof}

The construction of $f$ is the following.  

First solve the initial value problem for the wave equation:
$$\displaystyle
\left\{
\begin{array}{ll}
(\partial_s^2-\Delta)v=0 & \text{in $\Bbb R^3\times\,]0,\,\infty[$,}
\\
\\
\displaystyle
v(x,0)=0 & \text{in $\Bbb R^3$,}
\\
\\
\displaystyle
\partial_s v(x,0)=\Psi_B(x) & \text{in $\Bbb R^3$,}
\end{array}
\right.
\tag {2.1}
$$
where
$$\begin{array}{ll}
\displaystyle
\Psi_B(x)=(\eta-\vert x-p\vert)\chi_B(x), & x\in\Bbb R^3
\end{array}
$$
and the function $\chi_B(x)$ denotes the characteristic function of $B$.
Note that one can find the $v$ in the class
$$\displaystyle
C^2([0, \infty[, L^2(\Bbb R^3))\cap C^1([0,\,\infty[, H^1(\Bbb R^3))\cap C([0,\,\infty[, H^2(\Bbb R^3)).
$$
This is an application of the theory of $C_0$-semigroups \cite{Y}.
Note that the $v$ has an explicit form which is a consequence of {\it Kirchhoff's formula} \cite{ES}, that is
$$\displaystyle
v(x,s)=\frac{1}{4\pi s}
\,\int_{\partial B_s(x)}\,\Psi_B(y)dS_y,
$$
where $B_s(x)=\{y\in\Bbb R^3\,\vert\,\vert y-x\vert<s\}$.

Let $\tau>0$ and define
$$\begin{array}{lll}
\displaystyle
f_{B,\,T,\,\tau}(x,t)=\frac{1}{\sqrt{\tau}}\cdot\frac{\partial}{\partial\nu}v(x,\sqrt{\tau}(T-t)),
& x\in\partial\Omega, & 0<t<T.
\end{array}
$$
Note that the function $v_{\tau}(x,t), x\in\Bbb R^3, t\in\,[0,\,T]$ defined by
$$\displaystyle
v_{\tau}(x,t)=\frac{1}{\sqrt{\tau}}v(x,\sqrt{\tau}t),
$$
satisfies
$$\displaystyle
\left\{
\begin{array}{ll}
(\partial_t^2-\tau\Delta)v_{\tau}=0 & \text{in $\Bbb R^3\times\,]0,\,T[$,}
\\
\\
\displaystyle
v_{\tau}(x,0)=0 & \text{in $\Bbb R^3$,}
\\
\\
\displaystyle
\partial_t v_{\tau}(x,0)=\Psi_B(x) & \text{in $\Bbb R^3$.}
\end{array}
\right.
\tag {2.2}
$$
Thus one has the expression
$$\begin{array}{lll}
\displaystyle
f_{B,\,T,\,\tau}(x,t)=\frac{\partial}{\partial\nu}v_{\tau}(x,T-t),
& x\in\partial\Omega, & 0<t<T.
\end{array}
$$
This means that $f_{B,\,T,\,\tau}$ is the time-reversal mirror \cite{FWCM} on $\partial\Omega$
of the wave $v_{\tau}(x,t)$.

As can be seen in \cite{PW} it is not difficult to show that
$$\displaystyle
\Vert f_{B,\,T,\,\tau}(\,\cdot\,,t)\Vert_{H^{1/2}(\partial\Omega)}
\le C_{\Omega,B}\sqrt{(T-t)^2+\frac{3}{\tau}}
$$
and hence $f_{B,\,T,\,\tau}=O(1)$ in $L^2(0,T; H^{-1/2}(\partial\Omega))$ as $\tau\rightarrow\infty$.
Therefore from (1.3) we have, for $f=f_{B,\,T,\,\tau}$
$$\displaystyle
\Vert u_f(\,\cdot\,,T)\Vert_{L^2(\Omega\setminus\overline D)}=O(1).
\tag {2.3}
$$

Using $u_f$ with $f=f_{B,\,T,\,\tau}$, we define the indicator function
by the formula
$$\displaystyle
I_{\partial\Omega}(\tau;B,T)
=\int_{\partial\Omega}(w_f-w_f^*)\frac{\partial w_f^*}{\partial\nu}\,dS,
$$
where
$$\displaystyle
\left\{
\begin{array}{ll}
\displaystyle
w_f(x,\tau)=\int_0^T e^{-\tau t} u_f(x,t)dt, & x\in\Omega\setminus\overline D,
\\
\\
\displaystyle
w_f^*(x,\tau)=\int_0^T e^{-\tau t} v_{\tau}(x,T-t)\,dt, & x\in\Bbb R^3.
\end{array}
\right.
$$
Now we state the result of this paper.

\proclaim{\noindent Theorem 2.1.}
Let $\eta$ satisfy
$$\displaystyle
\eta+2R_D(p)>R_{\Omega}(p).
\tag {2.4}
$$

(i)  There exists a positive number $\tau_0$ such that $I_{\partial\Omega}(\tau;B,T)>0$
for all $\tau\ge\tau_0$ and
we have
$$\displaystyle
\lim_{\tau\rightarrow\infty}
\frac{1}{\sqrt{\tau}}\left(\log I_{\partial\Omega}(\tau;B,T)+2\tau T\right)
=2(\eta+R_D(p)).
\tag {2.5}
$$

(ii) We have
$$\displaystyle
\lim_{\tau\rightarrow\infty}
e^{-\sqrt{\tau}\,T}e^{2\tau T}I_{\partial\Omega}(\tau;B,T)
=
\left\{
\begin{array}{ll}
0 & \text{if $T>2(\eta+R_D(p))$,}
\\
\\
\displaystyle
\infty & \text{if $T<2(\eta+R_D(p))$.}
\end{array}
\right.
$$

\endproclaim

{\bf\noindent  Remark 2.2.}
The $\eta$ which is the radius of ball $B$ can not be arbitrary small since we have the constraint (2.4)
and it is really constraint if $2R_D(p)<R_{\Omega}(p)$. 
If $2R_D(p)\ge R_{\Omega}(p)$, then $\eta$ can be arbitrary small.
To cover the both cases, one can choose an arbitrary $\eta$ such that $\eta\ge R_{\Omega}(p)$.  
This is an advantage of using the initial data $\Psi_B$ in (2.1) with finitely extended support $\overline B$ not like $\{p\}$.

{\it\noindent Proof.}
Note that the assertion (ii) is a direct consequence from (i) since there is no restriction on $T$ in (i).
Let $w^*=w_f^*$.
One can write
$$\displaystyle
w^*(x,\tau)=e^{-\tau T}w_0(x,\tau),
$$
where
$$\displaystyle
w_0(x,\tau)=\int_0^T e^{\tau s}v_{\tau}(x,s)ds.
$$
From (2.2) we have
$$\displaystyle
(\Delta-\tau)w_0
+\frac{1}{\tau}\Psi_B(x)
=\frac{e^{\tau T}}{\tau}
(\partial_tv_{\tau}(x,T)-\tau v_{\tau}(x,T)).
$$
Thus, $w^*$ satisfies
$$\begin{array}{ll}
\displaystyle
(\Delta-\tau)w^*+\frac{1}{\tau}
(\sqrt{\tau} v(x,\sqrt{\tau}\,T)-\partial_sv(x,\sqrt{\tau}\,T))
=e^{-\tau T}F_0(x), & x\in\Bbb R^3
\end{array}
\tag {2.6}
$$
where
$$\displaystyle
F_0(x)=-\frac{\Psi_B(x)}{\tau}.
\tag {2.7}
$$
One the other hand, from (1.1) we see that $w=w_f$ satisfies
$$\displaystyle
\left\{
\begin{array}{ll}
(\Delta-\tau)w=e^{-\tau T}F(x) & \text{in $\Omega\setminus\overline D$,}
\\
\\
\displaystyle
\frac{\partial w}{\partial\nu}=\frac{\partial w^*}{\partial\nu} & \text{on $\partial\Omega$,}
\\
\\
\displaystyle
\frac{\partial w}{\partial\nu}=0 & \text{on $\partial D$,}
\end{array}
\right.
\tag {2.8}
$$
where
$$\begin{array}{ll}
\displaystyle
F(x)=u_f(x,T), & x\in\Omega\setminus\overline D.
\end{array}
\tag {2.9}
$$
Hereafter let $\tau$ satisfy
$$\displaystyle
\sqrt{\tau}\,T-\eta\ge R_{\Omega}(p),
$$
that is
$$\displaystyle
\tau\ge \left(\frac{\eta+R_{\Omega}(p)}{T}\right)^2.
$$
Then, we have $\Omega\subset B_{\sqrt{\tau}\,T-\eta}(p)$.
Since we have, for all $s>\eta$
$$\displaystyle
\text{supp}\,v(\,\cdot\,,s)\cup\text{supp}\,\partial_s v(\,\cdot\,,s)
\subset \Bbb R^3\setminus B_{s-\eta}(p),
\tag {2.10}
$$
from (2.6) we have
$$\begin{array}{ll}
\displaystyle
(\Delta-\tau)w^*=e^{-\tau T}F_0(x),
& \text{$x\in\Omega$}.
\end{array}
\tag {2.11}
$$
Note that (2.10) is a consequence of Kirchhoff's formula.

Set $R=w-w^*$.  Then from (2.8) and (2.11) we have
$$\displaystyle
\left\{
\begin{array}{ll}
(\Delta-\tau)R=e^{-\tau T}(F(x)-F_0(x)) & \text{in $\Omega\setminus\overline D$,}
\\
\\
\displaystyle
\frac{\partial R}{\partial\nu}=0 & \text{on $\partial\Omega$,}
\\
\\
\displaystyle
\frac{\partial R}{\partial\nu}=-\frac{\partial w^*}{\partial\nu} & \text{on $\partial D$.}
\end{array}
\right.
$$
Then, we obtain the following decomposition formula
which corresponds to (2.6) in \cite{PW}:
$$\displaystyle
I_{\partial\Omega}(\tau;B)
=J_h(\tau)+E_h(\tau)+{\cal R}_h(\tau),
\tag {2.12}
$$
where
$$
\displaystyle
\left\{
\begin{array}{l}
\displaystyle
J_h(\tau)=\int_D(\vert\nabla w^*\vert^2+\tau\vert w^*\vert^2)\,dx,\\
\\
\displaystyle
E_h(\tau)
=\int_{\Omega\setminus\overline D}(\vert\nabla R\vert^2+\tau\vert R\vert^2)\,dx
\end{array}
\right.
$$
and
$$\displaystyle
{\cal R}_h(\tau)
=e^{-\tau T}
\left\{\int_DF_0 w^*dx+\int_{\Omega\setminus\overline D}FRdx+
\int_{\Omega\setminus\overline D}
(F_0-F)w^*dx\right\}.
$$
And also similarly to \underline{Lemma 2.2} in \cite{PW} one can derive the estimate
$$\displaystyle
E_h(\tau)=O(\tau J_h(\tau)+e^{-2\tau T}).
\tag {2.13}
$$
From (2.6) and (2.7) we have the expression
$$\displaystyle
w^*=\frac{1}{\tau}\left(w_1^*+e^{-\tau T}w_R^*\right),
\tag {2.14}
$$
where
$$\begin{array}{ll}
\displaystyle
w_1^*(x,\tau)=\frac{1}{4\pi}\int_{\Bbb R^3}\frac{e^{-\sqrt{\tau}\vert x-y\vert}}{\vert x-y\vert}\,
(\sqrt{\tau}v(y,\sqrt{\tau}\,T)-\partial_s v(y,\sqrt{\tau}\,T))dy,
& x\in\Bbb R^3
\end{array}
\tag {2.15}
$$
and $w_R^*=w_R^*(x,\tau)$ satisfies
$$\begin{array}{ll}
\displaystyle
(\Delta-\tau)w_R^*+\Psi_B=0, & x\in\Bbb R^3.
\end{array}
$$
By integration by parts we have immediately, as $\tau\rightarrow\infty$,
$$
\displaystyle
\sqrt{\tau}\Vert w_R^*\Vert_{L^2(\Bbb R^3)}+\Vert\nabla w_R^*\Vert_{L^2(\Bbb R^3)}
=O(1).
\tag {2.16}
$$

Here we note that $v$ satisfies
$$\displaystyle
\text{supp}\,v(\,\cdot\,,s)\cup\text{supp}\,\partial_s v(\,\cdot\,,s)
\subset B_{s+\eta}(p)
\tag {2.17}
$$
for all $s\ge 0$, which is also a consequence of Kirchhoff's formula.
From (2.10) and (2.17) together with \underline{Propositions 3.3 and 3.4} in \cite{EV} we have
the explicit expression of the right-hand side on (2.15):
$$\displaystyle
w_1^*(x,\tau)
=\frac{1}{\tau}
({\cal H}_+(\sqrt{\tau};\sqrt{\tau}\,T,\eta)
+{\cal H}_{-}(\sqrt{\tau};\sqrt{\tau}\,T,\eta))
\,\frac{\sinh\sqrt{\tau}\vert x-p\vert}{\vert x-p\vert},
\tag {2.18}
$$
for all $x\in B_{\sqrt{\tau}\,T-\eta}(p)\setminus\{p\}$, where ${\cal H}_+(\sqrt{\tau};\sqrt{\tau}\,T,\eta)$ 
and ${\cal H}_{-}(\sqrt{\tau};\sqrt{\tau}\,T,\eta)$ are given by the following ${\cal H}_{+}(\tau;T,\eta)$ and ${\cal H}_{-}(\tau;T,\eta)$, respectively
in which $\tau$ and $T$ are replaced with
$\sqrt{\tau}$ and $\sqrt{\tau}\,T$, respectively:
$$\left\{
\begin{array}{l}
\displaystyle
{\cal H}_{+}(\tau;T,\eta)
=f_{\tau}(T)e^{-\tau T}-f_{\tau}(T+\eta)e^{-\tau (T+\eta)},
\\
\\
\displaystyle
{\cal H}_{-}(\tau;T,\eta)
=g_{\tau}(T-\eta)e^{-\tau (T-\eta)}-g_{\tau}(T)e^{-\tau T}.
\end{array}
\right.
$$
\underline{Here $f_{\tau}(\,\cdot\,)$ and $g_{\tau}(\,\cdot\,)$ are 
cubic polynomials given by}
$$\begin{array}{ll}
\displaystyle
f_{\tau}(\xi)
&
\displaystyle
=\frac{\tau}{6}\xi^3+\left\{1-\frac{\tau}{4}(\eta+2T)\right\}\xi^2+\left\{\frac{1}{2}\tau T(\eta+T)-(\eta+2T)+\frac{2}{\tau}\right\}\xi
\\
\\
\displaystyle
&
\displaystyle
\,\,\,
+\left\{\frac{1}{12}\tau(\eta-2T)(\eta+T)^2
+T(\eta+T)-\frac{\eta+2T}{\tau}+\frac{2}{\tau^2}\right\}
\end{array}
$$
and
$$\begin{array}{ll}
\displaystyle
g_{\tau}(\xi)
&
\displaystyle
=-\frac{\tau}{6}\xi^3-\left\{1+\frac{\tau}{4}(\eta-2T)\right\}\xi^2+\left\{\frac{1}{2}\tau T(\eta-T)-(\eta-2T)-\frac{2}{\tau}\right\}\xi\\
\\
\displaystyle
&
\displaystyle
\,\,\,
+\left\{\frac{1}{12}\tau(\eta+2T)(\eta-T)^2
+T(\eta-T)-\frac{\eta-2T}{\tau}-\frac{2}{\tau^2}\right\}.
\end{array}
$$
Note that, in the derivation of (2.18) \underline{we} fully made use of Kirchhoff's formula which gives the solution form
of (2.1) in the $(x,t)$-space.
Furthermore we have
$$\displaystyle
g_{\tau}(T-\eta)=\frac{1}{\tau}\left(\eta-\frac{2}{\tau}\right).
\tag {2.19}
$$
See (3.8) in \cite{EV}.
Using these, we have
$$\displaystyle
e^{\sqrt{\tau}(\sqrt{\tau}\,T-\eta)}(
{\cal H}_+(\sqrt{\tau};\sqrt{\tau}\,T,\eta)
+{\cal H}_{-}(\sqrt{\tau};\sqrt{\tau}\,T,\eta))
=\frac{1}{\sqrt{\tau}}
\left(\eta+O(\tau^{-1/2})\right).
$$
Then, applying \underline{Lemma 2.4} in \cite{EV} to the right-hand side on (2.18), we obtain, as
$\tau\rightarrow\infty$
$$\displaystyle
e^{2\tau T-2\sqrt{\tau}\,\eta}\left(\Vert w_1^*\Vert_{L^2(U)}^2+\Vert\nabla w_1^*\Vert_{L^2(U)}^2\right)=O(\tau^{\mu_1}e^{2\sqrt{\tau}\,R_U(p)})
\tag {2.20}
$$
and
$$\displaystyle
e^{2\tau T-2\sqrt{\tau}\,\eta}
\tau^{\mu_2}e^{-2\sqrt{\tau}\,R_U(p)}\,\Vert w_1^*\Vert_{L^2(U)}^2\ge C,
\tag {2.21}
$$
where $U$ is an arbitrary bounded open subset of $\Bbb R^3$; $p$ an arbitrary point in $\Bbb R^3$;
$\mu_1$ and $\mu_2$ are real numbers and $C$ a positive constant; it is assumed that
$\partial U$ is Lipschitz for (2.21).

From (2.14), (2.16) and (2.20) we obtain, as $\tau\rightarrow\infty$
$$\begin{array}{l}
\displaystyle
\,\,\,\,\,\,
e^{2\tau T-2\sqrt{\tau}\,\eta}(\Vert w^*\Vert_{L^2(U)}^2+\Vert\nabla w^*\Vert_{L^2(U)}^2)\\
\\
\displaystyle
=O(\tau^{-2}\tau^{\mu_1}e^{2\sqrt{\tau}\,R_U(p)})+O(\tau^{-2}e^{2\tau T-2\sqrt{\tau}\,\eta}e^{-2\tau T})\\
\\
\displaystyle
=O(\tau^{\mu_1-2}e^{2\sqrt{\tau}\,R_U(p)})(1+\tau^{-\mu_1}e^{-2\sqrt{\tau}\,\eta}e^{-2\sqrt{\tau}\,R_U(p)}))\\
\\
\displaystyle
=O(\tau^{\mu_1-2}e^{2\sqrt{\tau}\,R_U(p)}).
\end{array}
\tag{2.22}
$$
From (2.14), (2.16) and (2.21) we obtain
$$\begin{array}{l}
\displaystyle
\,\,\,\,\,\,
e^{2\tau T-2\sqrt{\tau}\,\eta}e^{-2\sqrt{\tau}\,R_U(p)}\Vert w^*\Vert_{L^2(U)}^2
\\
\\
\displaystyle
\ge \tau^{-2}\left(\frac{1}{2}e^{2\tau T-2\sqrt{\tau}\,\eta}e^{-2\sqrt{\tau}\,R_U(p)}\Vert w_1^*\Vert_{L^2(U)}^2
-e^{2\tau T-2\sqrt{\tau}\,\eta}e^{-2\sqrt{\tau}\,R_U(p)}e^{-2\tau T}\Vert w_R^*\Vert_{L^2(U)}^2
\right)\\\\
\displaystyle
\ge
\tau^{-2}\left(\frac{C}{2}\,\tau^{-\mu_2}
-C_1e^{-2\sqrt{\tau}\,\eta}e^{-2\sqrt{\tau}\,R_U(p)}\tau^{-1}\right)
\\
\\
\displaystyle
\ge C_2\tau^{-(2+\mu_2)},
\end{array}
\tag {2.23}
$$
where $C_1$, $C_2$ are positive constants being independent of all $\tau\ge\tau_0$ and $\tau_0$ is a large positive constant.

From (2.22) with $U=D$ we have
$$\displaystyle
e^{2\tau T-2\sqrt{\tau}\,\eta}e^{-2\sqrt{\tau}\,R_D(p)}J_h(\tau)
=O(\tau^{\mu_1}).
\tag {2.24}
$$
Then from (2.13) we have
$$\begin{array}{ll}
\displaystyle
e^{2\tau T-2\sqrt{\tau}\,\eta}e^{-2\sqrt{\tau}\,R_D(p)}E_h(\tau)
&
\displaystyle
=O(\tau^{\mu_1+1})
+O(e^{-2\sqrt{\tau}\,\eta}e^{-2\sqrt{\tau}\,R_D(p)})
\\
\\
\displaystyle
&
\displaystyle
=O(\tau^{\mu_1+1}).
\end{array}
\tag {2.25}
$$
Moreover, applying (2.3) to (2.9) and using (2.7), we have
$$\begin{array}{l}
\displaystyle
\,\,\,\,\,\,
\left\vert{\cal R}_h(\tau)\right\vert
\\
\\
\displaystyle
\le 
e^{-\tau T}
\left(\Vert F_0\Vert_{L^2(D)}\Vert w^*\Vert_{L^2(D)}+\Vert F\Vert_{L^2(\Omega\setminus\overline D)}
\Vert R\Vert_{L^2(\Omega\setminus\overline D)}
+\Vert F_0-F\Vert_{L^2(\Omega\setminus\overline D)}\Vert w^*\Vert_{L^2(\Omega\setminus\overline D)}\right)
\\
\\
\displaystyle
\le
Ce^{-\tau T}
\left(
\tau^{-1}\Vert w^*\Vert_{L^2(D)}+
\Vert R\Vert_{L^2(\Omega\setminus\overline D)}
+\Vert w^*\Vert_{L^2(\Omega\setminus\overline D)}\right).
\end{array}
$$
Then, from (2.24), (2.25) and (2.22) with $U=\Omega\setminus\overline D$
we obtain
$$
\begin{array}{l}
\displaystyle
\,\,\,\,\,\,
e^{\tau T-\sqrt{\tau}\,\eta}\left\vert{\cal R}_h(\tau)\right\vert
\\
\\
\displaystyle
=
Ce^{-\tau T}
\left(
O(\tau^{\frac{\mu_1}{2}-2}e^{\sqrt{\tau}\,R_D(p)})
+
O(\tau^{\frac{\mu_1}{2}-1})e^{\sqrt{\tau}\,R_{D}(p)}
+O(\tau^{\frac{\mu_1}{2}-1}e^{\sqrt{\tau}\,R_{\Omega\setminus\overline D}(p)})
\right)
\\
\\
\displaystyle
=
Ce^{-\tau T}
\left(
O(\tau^{\frac{\mu_1}{2}-1})e^{\sqrt{\tau}\,R_{D}(p)}
+O(\tau^{\frac{\mu_1}{2}-1}e^{\sqrt{\tau}\,R_{\Omega}(p)})
\right).
\end{array}
$$
This yields
$$
\begin{array}{l}
\displaystyle
\,\,\,\,\,\,
e^{2\tau T-2\sqrt{\tau}\,\eta}e^{-2\sqrt{\tau}\, R_D(p)}\left\vert{\cal R}_h(\tau)\right\vert
\\
\\
\displaystyle
=
O(\tau^{\frac{\mu_1}{2}-1}e^{-\sqrt{\tau}\,\eta}e^{-\sqrt{\tau}\,R_{D}(p)})
+O(\tau^{\frac{\mu_1}{2}-1}e^{-\sqrt{\tau}\,(\eta+2R_D(p)-R_{\Omega}(p)}).
\end{array}
\tag {2.26}
$$
Now we make use of the assumption (2.4).
Then, from (2.26) one concludes
$$\displaystyle
e^{2\tau T-2\sqrt{\tau}\,\eta}e^{-2\sqrt{\tau}\, R_D(p)}{\cal R}_h(\tau)
=O(\tau^{-\infty}).
\tag {2.27}
$$
Applying this together with (2.24) and (2.25) to (2.12), we obtain, as $\tau\rightarrow\infty$
$$\displaystyle
e^{2\tau T-2\sqrt{\tau}\,\eta}e^{-2\sqrt{\tau}\, R_D(p)}I_{\partial\Omega}(\tau;B,T)
=O(\tau^{\mu_1+1}).
\tag {2.28}
$$

From (2.23) with $U=D$ we have
$$\begin{array}{l}
\displaystyle
\,\,\,\,\,\,
e^{2\tau T-2\sqrt{\tau}\,\eta}e^{-2\sqrt{\tau}\, R_D(p)}\tau^{\mu_2}J_h(\tau)
\\
\\
\displaystyle
\ge e^{2\tau T-2\sqrt{\tau}\,\eta}e^{-2\sqrt{\tau}\, R_D(p)}\tau^{2+\mu_2}\Vert w^*\Vert_{L^2(D)}^2
\\
\\
\displaystyle
\ge C_2.
\end{array}
$$
Thus this together with (2.12) and (2.27) gives, for all $\tau\ge\tau_0$
for a sufficiently large $\tau_0>0$
$$\displaystyle
e^{2\tau T-2\sqrt{\tau}\,\eta}e^{-2\sqrt{\tau}\, R_D(p)}\tau^{\mu_2}I_{\partial\Omega}(\tau;B,T)
\ge \frac{C_2}{2}.
\tag {2.29}
$$

Now the assertion (i) follows from (2.28) and (2.29).

\noindent
$\Box$

{\bf\noindent Remark 2.3.}
Let $T=2(\eta+R_D(p))$.  This is the excluded case in Theorem 2.1 (ii).
It follows from (2.28) that, as $\tau\rightarrow\infty$
$$\displaystyle
e^{-\sqrt{\tau}\,T}e^{2\tau T}I_{\partial\Omega}(\tau;B,T)=O(\tau^{\mu_1+1}).
$$
The point is: this right-hand side is at most {\it algebraic}.

{\bf\noindent Remark 2.4.}
It follows from (2.28) and (2.29) that, as $\tau\rightarrow\infty$
$$\displaystyle
\frac{1}{\tau}\,\log I_{\partial\Omega}(\tau;B,T)=-2T+\frac{2}{\sqrt{\tau}}\,(\eta+R_D(p))+O(\frac{\log \tau}{\tau}).
\tag {2.30}
$$
This yields the convergence rate for (2.5):
$$\displaystyle
\frac{1}{\sqrt{\tau}}
\,(\log I_{\partial\Omega}(\tau;B,T)+2\tau T)
=2(\eta+R_D(p))+O(\frac{\log \tau}{\sqrt{\tau}}).
$$

{\bf\noindent Remark 2.5.}
In \cite{PW} instead of $f_{B,T,\tau}$ we have made use of the heat flux $f=f_{B,\,\tau}$ given by
$$\begin{array}{lll}
\displaystyle
f_{B,\,\tau}(x,t)=\frac{\partial}{\partial\nu}v_{\tau}(x,t), & x\in\partial\Omega, & 0<t<T,
\end{array}
$$
where $B$ is the same open ball as above, however, satisfies the constraint $\overline B\cap\overline\Omega=\emptyset$.
Note that there is no additional constraint on $B$ like (2.4).
Therein it is shown that, replacing $w_f^*$ in $I_{\partial\Omega}(\tau;B,T)$ with $w_f^0$ 
given by
$$\begin{array}{ll}
\displaystyle
w_f^0(x,\tau)=\int_0^T e^{-\tau t}v_{\tau}(x,t)dt, & x\in\Bbb R^3,
\end{array}
$$
one can extract $\text{dist}\,(D,B)$
by the formula
$$\displaystyle
\lim_{\tau\rightarrow\infty}\frac{1}{\sqrt{\tau}}\,\log I_{\partial\Omega}(\tau;B,T)
=-2\,\text{dist}\,(D,B).
$$
Comparing this with (2.30), one may think that the information $R_D(p)$ is {\it hidden} deeply more than $\text{dist}\,(D,B)$.

\section{An alternative proof of (2.18)}

In \underline{Remark 3.5} of \cite{EV} we have already pointed out that: together with (2.19) the formula
$$
\left\{
\begin{array}{l}
\displaystyle
f_{\tau}(T)=\frac{\tau}{12}\eta^3-\frac{\eta}{\tau}+\frac{2}{\tau^2},\\
\\
\displaystyle
f_{\tau}(T+\eta)=\frac{1}{\tau}\left(\eta+\frac{2}{\tau}\right),\\
\\
\displaystyle
g_{\tau}(T)=\frac{\tau}{12}\eta^3-\frac{\eta}{\tau}-\frac{2}{\tau^2}
\end{array}
\right.
$$
is valid.  Therefore we have
$$\begin{array}{l}
\displaystyle
\,\,\,\,\,\,
{\cal H}_+(\tau,T,\eta)
+{\cal H}_{-}(\tau;T,\eta)
\\
\\
\displaystyle
=\frac{4}{\tau^2}e^{-\tau T}
-\frac{1}{\tau}\left(\eta+\frac{2}{\tau}\right)e^{-\tau(T+\eta)}
+\frac{1}{\tau}\left(\eta-\frac{2}{\tau}\right)e^{-\tau(T-\eta)}.
\end{array}
\tag {3.1}
$$
We do not need this explicit formula (3.1) for the present purpose.
However, for future purpose which aims at extending the method presented in this paper to
other partial differential equations, including systems,
we present here an another proof which is based on the solution formula on (2.1) in $(\xi,t)$-space not $(x,t)$-space.
It is a combination of Parseval's identity and the residue theorem and completely different from the original proof
given in \cite{EV} which is a combination of Kirchhoff's formula and explicit computation formula of some volume potentials.

Let us formulate what we give a proof.

\proclaim{\noindent Proposition 3.1.}
Let $v$ be the solution of (2.1).
Let $T>\eta$.
We have
$$\begin{array}{l}
\,\,\,\,\,\,\displaystyle
\frac{\tau^2}{4\pi}
\int_{\Bbb R^3}
\frac{e^{-\tau\vert x-y\vert}}{\vert x-y\vert}
(\tau v(y,T)-\partial_tv(y,T))\,dy\\
\\
\displaystyle
=
\left({\cal H}_+(\tau,T,\eta)
+{\cal H}_{-}(\tau;T,\eta)\right)
\,\frac{\sinh\tau\vert x-p\vert}{\vert x-p\vert}
\end{array}
$$
for all $x\in B_{T-\eta}(p)\setminus\{p\}$,  where ${\cal H}_+(\tau,T,\eta)+{\cal H}_{-}(\tau;T,\eta)$ 
is given by the right-hand side on (3.1).
\endproclaim

{\it\noindent Proof.}
It is known that
$$\displaystyle
\hat{v}(\xi,T)
=\frac{\sin\vert\xi\vert T}{\vert\xi\vert}\,\hat{\Psi_B}(\xi)
$$
and
$$\displaystyle
\hat{(\partial_tv)}(\xi,T)
=\cos\vert\xi\vert T \,\hat{\Psi_B}(\xi),
$$
where we denote by $\hat{h}$ the Fourier transfom of a fuinction $h(x)$, that is
$$\displaystyle
\hat{h}(\xi)=\int_{\Bbb R^3}e^{-ix\cdot\xi}h(x)dx.
$$
We have also
$$
\displaystyle
\frac{1}{4\pi}\int_{\Bbb R^3}\frac{e^{-\tau\vert x\vert}}{\vert x\vert}e^{-ix\cdot\xi}dx
=\frac{1}{\vert\xi\vert^2+\tau^2}
$$
and thus
$$
\displaystyle
\frac{1}{4\pi}\int_{\Bbb R^3}\frac{e^{-\tau\vert x-y\vert}}{\vert x-y\vert}e^{-iy\cdot\xi}dy
=\frac{1}{\vert\xi\vert^2+\tau^2}e^{-ix\cdot\xi}.
$$
By Parseval's identity, we obtain
$$\begin{array}{l}
\,\,\,\,\,\,
\displaystyle
I_1(\tau)\equiv
\frac{\tau^2}{4\pi}
\int_{\Bbb R^3}
\frac{e^{-\tau\vert x-y\vert}}{\vert x-y\vert}
\tau v(y,T)\,dy\\
\\
\displaystyle
=\frac{\tau^2}{(2\pi)^3}\int_{\Bbb R^3}\tau\frac{\sin\vert\xi\vert T}{\vert\xi\vert}
\,\hat{\Psi_B}(\xi)\frac{1}{\vert\xi\vert^2+\tau^2}e^{ix\cdot\xi}d\xi
\end{array}
$$
and
$$\begin{array}{l}
\,\,\,\,\,\,
\displaystyle
I_2(\tau)\equiv
\frac{\tau^2}{4\pi}
\int_{\Bbb R^3}
\frac{e^{-\tau\vert x-y\vert}}{\vert x-y\vert}
\partial_tv(y,T)\,dy\\
\\
\displaystyle
=\frac{\tau^2}{(2\pi)^3}\int_{\Bbb R^3}\cos\vert\xi\vert T
\,\hat{\Psi_B}(\xi)\frac{1}{\vert\xi\vert^2+\tau^2}e^{ix\cdot\xi}d\xi.
\end{array}
$$
A change of variables yields
$$\begin{array}{ll}
\displaystyle
\hat{\Psi_B}(\xi)
&
\displaystyle
=\frac{4\pi\eta^3e^{-i\xi\cdot p}}{\vert\xi\vert}
\int_0^1r(1-r)\sin\eta\vert\xi\vert r\,dr.
\end{array}
$$
Here we have
$$
\displaystyle
\int_0^1r(1-r)\sin\eta\vert\xi\vert r\,dr
=\frac{1}{(\eta\vert\xi\vert)^2}
\left\{\frac{2}{\eta\vert\xi\vert}\,(1-\cos\eta\vert\xi\vert)-\sin\eta\vert\xi\vert\right\}.
$$
Thus we obtain
$$\begin{array}{ll}
\displaystyle
\hat{\Psi_B}(\xi)
&
\displaystyle
=\frac{4\pi\eta^3e^{-i\xi\cdot p}}{\vert\xi\vert}\cdot\frac{1}{(\eta\vert\xi\vert)^2}
\left\{\frac{2}{\eta\vert\xi\vert}\,(1-\cos\eta\vert\xi\vert)-\sin\eta\vert\xi\vert\right\}\\
\\
\displaystyle
&
\displaystyle
=4\pi\eta e^{-i\xi\cdot p}g(\vert\xi\vert;\eta),
\end{array}
$$
where
$$\displaystyle
g(z;\eta)
=\,\frac{1}{z^3}
\left\{\frac{2}{\eta z}\,(1-\cos\eta z)-\sin\eta z\right\}.
$$
Note that $g(z;\eta)$ is an even function and that $z=0$ is a removable singularity of $g(z;\eta)$.

Therefore we obtain the expression
$$\begin{array}{ll}
\displaystyle
I_1(\tau)
&
\displaystyle
=\frac{\tau^2}{(2\pi)^3}\int_{\Bbb R^3}\tau\frac{\sin\vert\xi\vert T}{\vert\xi\vert}
\cdot 4\pi\eta e^{-i\xi\cdot p}g(\vert\xi\vert;\eta)\frac{1}{\vert\xi\vert^2+\tau^2}e^{ix\cdot\xi}d\xi
\\
\\
\displaystyle
&
\displaystyle
=\frac{\tau^3\eta}{2\pi^2}\int_{\Bbb R^3}\frac{\sin\vert\xi\vert T}{\vert\xi\vert}
\cdot
\frac{g(\vert\xi\vert;\eta)}{\vert\xi\vert^2+\tau^2}e^{i(x-p)\cdot\xi}d\xi\\
\\
\displaystyle
&
\displaystyle
=\frac{\tau^2\eta}{2\pi^2}\int_{\Bbb R^3}
\frac{\sin\vert\xi\vert T}{\vert\xi\vert}
\frac{g(\vert\xi\vert;\eta)}{\vert\xi\vert^2+\tau^2}e^{i\vert x-p\vert\xi_3}d\xi
\\
\\
\displaystyle
&
\displaystyle
=\frac{\tau^2\eta}{2\pi^2}
\int_0^{\infty}r^2 dr
\int_0^{2\pi}d\theta\int_0^{\pi}\sin\varphi d\varphi
\frac{\sin r T}{r}
\frac{g(r;\eta)}{r^2+\tau^2}e^{i\vert x-p\vert r\cos\varphi}\\
\\
\displaystyle
&
\displaystyle
=\frac{\tau^2\eta}{\pi}
\int_0^{\infty}\frac{r^2}{r^2+\tau^2}\frac{\sin r T}{r}g(r;\eta)
dr\int_0^{\pi}e^{i\vert x-p\vert r\cos\varphi}\sin\varphi d\varphi
\\
\\
\displaystyle
&
\displaystyle
=\frac{i\tau^2\eta}{\pi\vert x-p\vert}\int_0^{\infty}\frac{r}{r^2+\tau^2}
\frac{\sin r T}{r}g(r;\eta)
dr\int_0^{\pi}
(e^{i\vert x-p\vert r\cos\varphi})'d\varphi
\\
\\
\displaystyle
&
\displaystyle
=\frac{2\tau^2\eta}{\pi\vert x-p\vert}\int_0^{\infty}\frac{r}{r^2+\tau^2}
\frac{\sin r T}{r}g(r;\eta)
\sin\vert x-p\vert r \,dr
\\
\\
\displaystyle
&
\displaystyle
=\frac{\tau^2\eta}{\pi\vert x-p\vert}\int_{-\infty}^{\infty}\frac{\sin rT}{r^2+\tau^2}g(r;\eta)
\sin\vert x-p\vert r \,dr.
\end{array}
$$
Similarly, we have
$$
\displaystyle
I_2(\tau)
=\frac{\tau^2\eta}{\pi\vert x-p\vert}\int_{-\infty}^{\infty}\frac{r\cos rT}{r^2+\tau^2}
g(r;\eta)
\sin\vert x-p\vert r \,dr.
$$

Since $g(-r;\eta)=g(r;\eta)$, one can rewrite
$$
\left\{
\begin{array}{l}
\displaystyle
\int_{-\infty}^{\infty}\frac{\sin rT}{r^2+\tau^2}
g(r;\eta)
\sin\vert x-p\vert r\,dr
=\frac{1}{i}\int_{-\infty}^{\infty}\frac{e^{irT}}{r^2+\tau^2}
g(r;\eta)
\sin\vert x-p\vert r\,dr,\\
\\
\displaystyle
\int_{-\infty}^{\infty}\frac{r\cos rT}{r^2+\tau^2}
g(r;\eta)
\sin\vert x-p\vert r\,dr
=\int_{-\infty}^{\infty}\frac{re^{irT}}{r^2+\tau^2}
g(r;\eta)\sin\vert x-p\vert r\,dr.
\end{array}
\right.
$$
Therefore we have
$$\left\{
\begin{array}{l}
\displaystyle
I_1(\tau)
=-i\frac{\tau^3\eta}{\pi\vert x-p\vert}\int_{-\infty}^{\infty}\frac{e^{irT}}{r^2+\tau^2}
g(r;\eta)
\sin\vert x-p\vert r\,dr,
\\
\\
\displaystyle
I_2(\tau)
=\frac{\tau^2\eta}{\pi\vert x-p\vert}
\int_{-\infty}^{\infty}\frac{re^{irT}}{r^2+\tau^2}
g(r;\eta)\sin\vert x-p\vert r\,dr.
\end{array}
\right.
$$

Let $z=Re^{i\theta}$ with $0<\theta<\pi$, $R>0$.
We have, as $R\rightarrow\infty$
$$
\left\{
\begin{array}{l}
\displaystyle
\vert g(z;\eta)\vert=O(R^{-3}e^{\eta R\sin\theta}),\\
\\
\displaystyle
\vert\sin\vert x-p\vert z\vert=O(e^{\vert x-p\vert R\sin\theta}),\\
\\
\displaystyle
\vert e^{izT}\vert=O(e^{-TR\sin\theta}).
\end{array}
\right.
$$
Since we have $\displaystyle T>\eta+\vert x-p\vert$,
one can apply the residue theorem by taking a standard closed contour $-R\longrightarrow R\longrightarrow Re^{i\theta}\,(0\le\theta\le\pi)\longrightarrow -R$
with $R>\tau$ and letting $R\rightarrow\infty$, one gets
$$\begin{array}{ll}
\displaystyle
\int_{-\infty}^{\infty}\frac{e^{irT}}{r^2+\tau^2}
g(r;\eta)
\sin\vert x-p\vert r\,dr
&
\displaystyle
=2\pi i\,\text{Res}_{z=i\tau}\left(\frac{e^{izT}}{z^2+\tau^2}
g(z;\eta)
\sin\vert x-p\vert z\right)
\\
\\
\displaystyle
&
\displaystyle
=2\pi i
\left(\frac{e^{-\tau T}}{2i\tau}
g(i\tau;\eta)
\sin\vert x-p\vert i\tau\right)
\\
\\
\displaystyle
&
\displaystyle
=-\frac{\pi}{i\tau}e^{-\tau T}g(i\tau;\eta)\sinh \tau\vert x-p\vert.
\end{array}
$$
Similarly, we have
$$\begin{array}{ll}
\displaystyle
\int_{-\infty}^{\infty}\frac{re^{irT}}{r^2+\tau^2}
g(r;\eta)\sin\vert x-p\vert r\,dr
&
\displaystyle
=2\pi i\frac{i\tau e^{-\tau T}}{2i\tau}
g(i\tau;\eta)\sin\vert x-p\vert i\tau
\\
\\
\displaystyle
&
\displaystyle
=-\pi e^{-\tau T}
g(i\tau;\eta)\sinh\tau\vert x-p\vert.
\end{array}
$$
From these one gets
$$
\displaystyle
I_1(\tau)
=\tau^2\eta
e^{-\tau T}
g(i\tau;\eta)\frac{\sinh\tau\vert x-p\vert}{\vert x-p\vert}
$$
and $I_2(\tau)=-I_1(\tau)$.

Therefore we obtain
$$\displaystyle
I(\tau)=2I_1(\tau)=
2\tau^2\eta e^{-\tau T}g(i\tau;\eta)\frac{\sinh\tau\vert x-p\vert}{\vert x-p\vert}.
$$
Here we have
$$\begin{array}{ll}
\displaystyle
g(i\tau;\eta)
&
\displaystyle
=\frac{1}{(i\tau)^3}
\left\{\frac{2}{\eta i\tau}\,(1-\cos i\eta\tau)-\sin i\eta\tau\right\}
\\
\\
\displaystyle
&
\displaystyle
=\frac{2}{\eta\tau^4}\left(1-\frac{e^{\tau\eta}+e^{-\tau\eta}}{2}\right)-\frac{1}{(i\tau)^3}\frac{e^{-\tau\eta}-e^{\tau\eta}}{2i}\\
\\
\displaystyle
&
\displaystyle
=\frac{2}{\eta\tau^4}\left(1-\frac{e^{\tau\eta}+e^{-\tau\eta}}{2}\right)-\frac{1}{\tau^3}\frac{e^{-\tau\eta}-e^{\tau\eta}}{2}\\
\\
\displaystyle
&
\displaystyle
=\frac{2}{\eta\tau^4}+\left(\frac{1}{2\tau^3}-\frac{1}{\eta\tau^4}\right)e^{\tau\eta}
-\left(\frac{1}{2\tau^3}+\frac{1}{\eta\tau^4}\right)e^{-\tau\eta}.
\end{array}
$$
This completes the proof.

\noindent
$\Box$

$$\quad$$

\centerline{{\bf Acknowledgment}}

The author was partially supported by Grant-in-Aid for
Scientific Research (C)(No. 17K05331) of Japan  Society for
the Promotion of Science.

$$\quad$$

\vskip1cm
\noindent
e-mail address

ikehata@hiroshima-u.ac.jp


\begin{thebibliography}{99}






\bibitem{BEL}  Belishev, M. I.,
             On an approach to multidimensional inverse problems for the wave equation,
             Dokl. Akad. Nauk SSSR, {\bf 297}(1987), 524-527.




\bibitem{B2} Belishev, M. I.,
             How to see waves under the Earth surface (the BC-method for geophysicists),
             Ill-Posed and Inverse Problems, pp. 67-84,
             Kabanikhin, S. I. and Romanov, V. G. (Eds), VSP, Utrecht, 2002.






\bibitem{BKLS} Bingham, K., Kurylev, Y., Lassas, M. and Siltanen, S.,
              Iterative time-reversal control for inverse problems, 
              Inverse Problems and Imaging, {\bf 2}(2008), 63-81.









\bibitem{dBK} de Buhan, M. and Kray, M., 
               A new approach to solve the inverse scattering problem for waves: 
               combining the TRAC and the Adaptive Inversion methods,
               Inverse Problems, {\bf 29}(2013), 085009.
               
               
               











\bibitem{DKL} Dahl, M. F., Kirpichnikova, A. and Lassas, M.,
              Focusing waves in unknown media by modified time reversal iteration,
              SIAM J. Control. Optim., {\bf 48}(2009), 839-858.









\bibitem{DL}  Dautray, R. and Lions, J-L., Mathematical analysis and numerical methods for
              sciences and technology.Vol. {\bf 5}.
              Evolution problems. I, Springer-Verlag, Berlin, 1992.










\bibitem{ES}  Egorov, Yu. V. and Shubin, M. A., 
              Foundations of the classical theory of partial differential equations,
              Springer, 1998.



\bibitem{FWCM} Fink, M.,
               Time reversal of ultrasonic fields-Part I: Basic principles,
               IEEE Trans. Ultrason., Ferroelec., Freq. Contr., {\bf 39}(1992), No.5, pp.555-566.










\bibitem{E001} Ikehata, M.,
               \newblock Enclosing a polygonal cavity in a two-dimensional bounded domain from Cauchy data,
               Inverse Problems, {\bf 15}(1999), 1231-1241.




\bibitem{E00} Ikehata, M.,
              \newblock Reconstruction of the support function for inclusion from boundary measurements,
              J. Inverse Ill-Posed Probl., {\bf 8}(2000), No. 4, 367-378.









\bibitem{I4} Ikehata, M.,
             \newblock Extracting discontinuity in a heat conductive body. One-space dimensional case,
             \newblock Applicable Analysis, {\bf 86}(2007), no. 8, 963-1005.








\bibitem{IW00} Ikehata, M.,
               The enclosure method for inverse obstacle scattering problems with dynamical data over a
               finite time interval, Inverse Problems, {\bf 26}(2010) 055010(20pp).





\bibitem{IFR} Ikehata, M.,
              The framework of the enclosure method with dynamical data and its applications,
              Inverse Problems, {\bf 27}(2011) 065005(16pp).






\bibitem{IEO2} Ikehata, M.,
              The enclosure method for inverse obstacle scattering problems with dynamical data over a finite time
              interval: II. Obstacles with a dissipative boundary or finite refractive index and back-scattering data,
              Inverse Problems, {\bf 28}(2012) 045010 (29pp).






\bibitem{IEO3} Ikehata, M.,
               The enclosure method for inverse obstacle scattering problems with dynamical data over a finite time interval:
               III. Sound-soft obstacle and bistatic data, Inverse Problems, {\bf 29}(2013) 085013 (35pp).













\bibitem{Iwall} Ikehata, M., 
                On finding an obstacle embedded in the rough background medium via the enclosure method in the time domain,
                Inverse Problems, {\bf 31}(2015) 085011(21pp).











\bibitem{IMax}  Ikehata, M.,
             The enclosure method for inverse obstacle scattering using a single electromagnetic
             wave in time domain,  Inverse Problems and Imaging, {\bf 10}(2016), 131-163.





\bibitem{IMax2} Ikehata, M.,
             On finding an obstacle with the Leontovich boundary condition via the time domain
             enclosure method, Inverse Problems and Imaging, {\bf 11}(2017), 99-123.







\bibitem{EIV} Ikehata, M., 
              The enclosure method for inverse obstacle scattering over a finite time interval: IV.
              Extraction from a single point on the graph of the response operator,
              J. Inverse Ill-Posed Probl., {\bf 25}(2017), No. 6, 747-761.
              




\bibitem{Ithermo} Ikehata, M.,
             On finding a cavity in a thermoelastic body using a single displacement measurement over a finite time interval on the surface of the body,
             J. Inverse Ill-Posed Probl., {\bf 26}(2018), No. 3, 369-394.

 
  


\bibitem{EV} Ikehata, M., 
             The enclosure method for inverse obstacle scattering over a finite time interval: V. Using
             time-reversal invariance, J. Inverse Ill-posed Probl., {\bf 27}(2019), No. 1, 133-149.
   
 
         

\bibitem{PW} Ikehata, M.,
              Prescribing a heat flux coming from a wave equation,
              J. Inverse Ill-Posed Probl., {\bf 27}(2019), No. 5, 731-744.
              
        
  

\bibitem{II}  Ikehata, M. and Itou, H., On reconstruction of a cavity in a linearized viscoelastic body from infinitely many transient
              boundary data, Inverse Problems, {\bf 28}(2012) 125003 (19pp).





\bibitem{IK1}  Ikehata, M. and Kawashita, M., The enclosure method for the heat equation,
               Inverse Problems, {\bf 25}(2009) 075005(10pp).
               
               
               



\bibitem{IK2}  Ikehata, M. and Kawashita, M., On the reconstruction of inclusions in a heat conductive body
               from dynamical boundary data over a finite time interval, Inverse Problems,
               {\bf 26}(2010) 095004(15pp).
    


\bibitem{IK3}  Ikehata, M. and Kawashita, M.,
               An inverse problem for a three-dimensional heat equation 
               in thermal imaging and the enclosure method, 
               Inverse Problems and Imaging, {\bf 8}(2014), 1073-1116.













\bibitem{O} Oksanen, L., 
             Solving an inverse obstacle problem for the wave equation by using the boundary control method,
             Inverse Problems, {\bf 29}(2013) 035004.








\bibitem{Y}   Yosida, K., Functional Analysis, Third Edition, Springer, New York, 1971.











\end{thebibliography}
\end{document}